\documentclass[10pt,leqno]{article}

\usepackage{amsmath,amsthm,amsfonts,latexsym,amscd,amssymb}
\numberwithin{equation}{section}
\input xypic
\def\qed{{\hbadness=10000\hfill\ \vbox{\hrule height.09ex
   \hbox{\vrule width.09ex height1.55ex depth.2ex \kern1.8ex
   \vrule width.09ex height1.55ex depth.2ex}\hrule height.09ex}\break
   \bigskip}}
\newtheorem{theorem}{Theorem}[section]

\theoremstyle{definition}

\theoremstyle{remark}

\begin{document}
\date{}

\linespread{1}\title{\textbf{\LARGE On a K\"ahlerian space-time manifold}} 
\author{\textbf{$^{*}$B. B. Chaturvedi  and $^{**}$Pankaj Pandey}\\
\normalsize{Department of Pure and Applied Mathematics}\\
 \normalsize{Guru Ghasidas Vishwavidyalaya Bilaspur (C.G.), India}\\
\normalsize{Email: $^{*}$brajbhushan25@gmail.com }\\
\normalsize{\:\:\:$^{**}$pankajpandey.mathematics@gmail.com}}
\maketitle

\noindent\textbf{Abstract} In this paper, the theory of space-time in 4-dimensional K\"ahler manifold has been studied. We have also discussed the Einstein equation with cosmological constant in perfect fluid K\"ahlerian space-time manifold. By taking conformally flat perfect fluid K\"ahler space-time manifold, we have obtained interesting results related to sectional curvatures. In last two sections we have studied weakly symmetric and weakly Ricci symmetric perfect fluid K\"ahler space-time manifolds.\\

 Einstein equation. conformally flat manifold. weakly symmetric manifold. weakly Ricci symmetric manifold\\


\section{Introduction}
The study of space time is associated with 4-dimensional semi-Riemannian manifolds equipped with Lorentz metric $g$ having signature (-,+,+,+). B. O'Neill \cite{nel83} discussed  the application of semi-Riemannian geometry in the theory of relativity. The curvature structure of the space-time is studied by V. R. Kaigorodov \cite{kai14} in 1983. After then, these ideas of general relativity of space-time are extended by A. K. Roychaudhury, S. Banerji and A. Banerji \cite{roy92}, M. C. Chaki and S. Roy \cite{chk96}, A. A. Shaikh, Dae Won Yoon and S. K. Hui \cite{shk09}, U. C. De and G. C. Ghosh \cite{ucd04} and many other Differential Geometers and Physicist.

In 2004 U. C. De and G. C. Ghosh \cite{ucd04} considered the weakly Ricci symmetric space-time manifold and obtained some results related to it. The weakly symmetric and weakly Ricci symmetric manifolds are introduced by L. Tamassy and T. Q. Binh \cite{tqb93, ltq42}. Also, M. Prvanovic \cite{mpr95} and U. C. De and S. Bandyopdhayay \cite{dsb99} explained it with examples. L. Tamassy, U. C. De and T. Q. Binh \cite{ltu00} also found some interesting results related to weakly symmetric and weakly Ricci symmetric K\"ahler manifolds.

An $n$-dimensional Riemannian manifold is said to be a weakly symmetric if the curvature tensor $R$ of type (0,4) of the manifold satisfies
 \begin{equation}
\begin{split}
(\nabla_XR)(Y,Z,U,V)=& A(X)R(Y,Z,U,V)\\&+B(Y)R(X,Z,U,V)+C(Z)R(Y,X,U,V)\\&+D(U)R(Y,Z,X,V)+E(V)R(Y,Z,U,X),
\end{split}
\end{equation}
and if the Ricci tensor $S$ of the manifold satisfies
\begin{equation}
(\nabla_XS)(Y,Z)=A(X)S(Y,Z)+B(Y)S(X,Z)+C(Z)S(Y,X),
\end{equation} 
then manifold is called weakly Ricci symmetric manifold. Where ${A,B,C,D,E}$ are simultaneously non-vanishing 1-forms and ${X,Y,Z,U,V}$ are vector fields.

In 1995, Prvanovic \cite{mpr95} proved that if the manifold be a weakly symmetric manifold satisfying (1.1) then $B=C=D=E$. In this paper we have taken $B=C=D=E=\omega$, therefore, the equations (1.1) and (1.2) can be written as
\begin{equation}
\begin{split}
(\nabla_XR)(Y,Z,U,V)=& A(X)R(Y,Z,U,V)\\&+\omega(Y)R(X,Z,U,V)+\omega(Z)R(Y,X,U,V)\\&+\omega(U)R(Y,Z,X,V)+\omega(V)R(Y,Z,U,X),
\end{split}
\end{equation}
and
\begin{equation}
(\nabla_XS)(Y,Z)=A(X)S(Y,Z)+\omega(Y)S(X,Z)+\omega(Z)S(Y,X),
\end{equation}
where $g(X,\rho)=\omega(X)$ and $g(X,\alpha)=A(X)$.

Recently, U. C. De and A. De \cite{dde12}, A. De, C. \"Ozg\"ur and U. C. De \cite{adc12} extended the theory of space-time by studying the application of almost pseudo-conformally symmetric Ricci-recurrent manifolds and conformally flat almost pseudo Ricci-symmetric manifolds respectively. In 2012 and 2014, S. Mallick and U. C. De \cite{uds14, smu14} gave examples of weakly symmetric and conformally flat weakly Ricci symmetric space-time manifolds and proved that a conformally flat weakly Ricci symmetric space-time with non-zero scalar curvature is the Robertson-Walker space-time and the vorticity and the shear vanish. Also, they have shown that a weakly symmetric perfect fluid space-time with cyclic parallel Ricci tensor can not admit heat flux. In 2014, U. C. De and Ljubica Velimirovic \cite{cdl14} explained the space-time with semi-symmetric energy-momentum tensor and proved that a space-time manifold with semi-symmetric energy momentum tensor is Ricci semi-symmetric.

In the consequence of this type of development we have studied the general relativity of space-time in 4-dimensional K\"ahler manifold and this type of manifold is called the K\"ahlerian space-time manifold.

A 4-dimensional space-time manifold is said to be a K\"ahlerian space-time manifold if the following conditions hold:
\begin{equation}
F^2(X)=-X,
\end{equation}
\begin{equation}
g(\overline X, \overline Y)=g(X,Y),
\end{equation}  
\begin{equation}
(\nabla_XF)(Y)=0,
\end{equation}
where, $F$ is a tensor field of type (1,1) such that ${F(X)=\overline X}$, ${g}$ is a Riemannian metric and ${\nabla}$ is a Levi-Civita connection.

We know that in a K\"ahler manifold the Ricci tensor S satisfies
\begin{equation}
S(\overline X, \overline Y)=S(X,Y).
\end{equation}

\section{Perfect fluid K\"ahler space-time manifold}

We know that the Einstein equation with cosmological constant for the perfect fluid space-time is given by
\begin{equation}
S(X,Y)-\frac{r}{2}g(X,Y)+\lambda g(X,Y)=k[(\sigma+p)\omega(X)\omega(Y)+pg(X,Y)],
\end{equation} 
where $k$ is the gravitational constant, $\sigma$ is energy density, $p$ is isotropic pressure of the fluid and $\omega$ is 1-form defined by $\omega(X)=g(X,\rho)$ for time-like vector field $\rho$. The time-like vector field $\rho$ is called velocity of the fluid and satisfies $g(\rho,\rho)=-1$. Also, the energy density $\sigma$ and the pressure $p$ can be described in the sense that if $\sigma$ vanishes then content matter of the fluid is not pure and if the pressure $p$ vanishes then the fluid is dust.

Now replacing X and Y by $\overline X$ and $\overline Y$ respectively in (2.1) and using (1.6) and (1.8), we get
\begin{equation}
S(X,Y)-\frac{r}{2}g(X,Y)+\lambda g(X,Y)=k[(\sigma+p)\omega(\overline X)\omega(\overline Y)+pg(X,Y)].
\end{equation} 
Subtracting (2.1) from (2.2), we have
\begin{equation}
k(\sigma+p)[\omega(\overline X)\omega(\overline Y)-\omega(X)\omega(Y)]=0.
\end{equation} 
Putting $Y=\rho$ in (2.3), we obtained
\begin{equation}
k(\sigma+p)\omega(X)=0,
\end{equation} 
since $k\neq 0$ and $\omega(X)\neq 0$, we have
\begin{equation}
\sigma + p=0,
\end{equation}
which shows the fluid behaves like a cosmological constant. Also from (2.5), we have $\sigma=-p$ which represents a rapid expansion of space-time in cosmology and known as Inflation.

Now using (2.5), equation (2.1) gives 
\begin{equation}
S(X,Y)=(\frac{r}{2}-\lambda + k.p)g(X,Y).
\end{equation} 
Putting $X=Y=e_i, 1\leq i\leq 4$ in (2.6) and taking summation over $i$, we can easily obtained
\begin{equation}
\lambda-k.p=\frac{r}{4}.
\end{equation} 
From (2.6) and (2.7), we get
\begin{equation}
S(X,Y)=\frac{r}{4}g(X,Y).
\end{equation}
Hence from above discussion, we conclude the following:
\begin{theorem}
If M be a perfect fluid K\"ahler space-time manifold satisfying the Einstein equation with cosmological constant then\\
(i) the Einstein equation is independent of the isotropic pressure $p$, energy density $\sigma$, cosmological constant $\lambda$  and the gravitational constant $k$.\\
(ii) the space-time manifold is an Einstein manifold.
\end{theorem}
Since an Einstein manifold is the manifold of constant scalar curvature $r$, therefore, equation (2.7) implies the pressure $p$ is constant and hence from (2.5) we have the energy density $\sigma$ is constant.

It is well known \cite{nel83} that the Energy equation for the perfect fluid is given by
 \begin{equation}
\rho.\sigma=-(\sigma+p)div\rho.
\end{equation} 
Using (2.5) in (2.9), we get
\begin{equation}
\rho.\sigma=0.
\end{equation}
Above equation implies $\sigma=0$ as $\rho\neq 0$. Because if $\rho= 0$ then we have contradiction that $g(\rho,\rho)=-1$. But then the equation (2.5) gives $p=0$ and hence the energy momentum tensor $T(X,Y)=(\sigma+p)\omega(X)\omega(Y)+pg(X,Y)$ vanishes.\\
Thus we conclude:
\begin{theorem}
If M be a perfect fluid K\"ahler space-time manifold satisfying the Einstein equation with cosmological constant then\\
(i) the energy density $\sigma$ and the isotropic pressure $p$ vanish $i.e.$ the content matter of the fluid is not pure and the perfect fluid is dust.\\
(ii) the energy momentum tensor $T(X,Y)$ vanishes $i.e.$ the space-time is vaccum.
\end{theorem}

 Since the velocity vector field $\rho$ is constant, we have 
\begin{equation}
div\rho=0 ~~~ and ~~~\nabla_\rho\rho=0.
\end{equation}
Hence from equations (2.11), we can state:
\begin{theorem}
If M be a perfect fluid K\"ahler space-time manifold satisfying Einstein equation with cosmological constant then the expansion scalar and the acceleration vector vanish.
\end{theorem}

\section{Conformally flat perfect fluid K\"ahler space-time manifold}

The Weyl conformal curvature tensor $C$ on an n-dimensional manifold $M$ is defined by
\begin{equation}
\begin{split}
C(X,Y,Z,T)=&R(X,Y,Z,T)-\frac{1}{(n-2)}[S(Y,Z)g(X,T)\\&-S(X,Z)g(Y,T)+S(X,T)g(Y,Z)\\&-S(Y,T)g(X,Z)]\\&+\frac{r}{(n-1)(n-2)}[g(Y,Z)g(X,T)-g(X,Z)g(Y,T)].
\end{split}
\end{equation}
If the manifold be conformally flat, then above equation implies
\begin{equation}
\begin{split}
R(X,Y,Z,T)=&\frac{1}{2}[S(Y,Z)g(X,T)-S(X,Z)g(Y,T)\\&+S(X,T)g(Y,Z)-S(Y,T)g(X,Z)]\\&-\frac{r}{6}[g(Y,Z)g(X,T)-g(X,Z)g(Y,T)].
\end{split}
\end{equation}
Putting the value of Ricci tensor S from (2.8) in (3.2), we have
\begin{equation}
R(X,Y,Z,T)=\frac{r}{12}[g(Y,Z)g(X,T)-g(X,Z)g(Y,T)].
\end{equation}
If the 3-dimensional distribution of the manifold orthogonal to $\rho$ is denoted by $\rho^\bot$  then by putting $Z=Y$ and $T=X$ in (3.3), we can write
\begin{equation}
R(X,Y,Y,X)=\frac{r}{12}[g(X,X)g(Y,Y)-g(X,Y)g(X,Y)],
\end{equation}
where X, Y $\in$ $\rho^\bot$.\\
Putting $Y=\rho$ in (3.4), we have
\begin{equation}
R(X,\rho,\rho,X)=-\frac{r}{12}g(X,X),
\end{equation}
where X $\in$ $\rho^\bot$.\\
If we denote the sectional curvatures determined by $X,Y$ and $X,\rho$ by $K(X,Y)$ and $K(X,\rho)$ respectively then from (3.4) and (3.5), we have
\begin{equation}
\begin{split}
K(X,Y)=&\frac{R(X,Y,Y,X)}{g(X,X)g(Y,Y)-g(X,Y)g(X,Y)}\\&=\frac{r}{12}.
\end{split}
\end{equation}
and
\begin{equation}
\begin{split}
K(X,\rho)=&\frac{R(X,\rho,\rho,X)}{g(X,X)g(\rho,\rho)-g(X,\rho)g(X,\rho)}\\&=\frac{r}{12}.
\end{split}
\end{equation} 
Thus we have:
\begin{theorem}
If M be a conformally flat perfect fluid K\"ahler space-time manifold satisfying the Einstein equation with cosmological constant then the sectional curvature determine by $X,Y$ and $X,\rho$ are same and equal to $\frac{r}{12}$  $i.e.$ $K(X,Y)=K(X,\rho)=\frac{r}{12}$. 
\end{theorem}
Karcher \cite{kar82} has defined that a Lorentzian manifold is said to be infinitesimally spatially isotropic relative to the velocity vector field $\rho$ if the Riemannian curvature tensor $R$ satisfies
 \begin{equation}
R(X,Y,Z,T)=a.[g(Y,Z)g(X,T)-g(X,Z)g(Y,T)],
\end{equation}
and
\begin{equation}
R(X,\rho,\rho,Y)=b.g(X,Y),
\end{equation} 
 where $a$ and $b$ are real valued functions and $X,Y,Z,T$ $\in$ $\rho^\bot$.
Putting $Y=Z=\rho$ in (3.3), we get
\begin{equation}
R(X,\rho,\rho,T)=-\frac{r}{12}g(X,T).
\end{equation}
Hence from (3.3) and (3.10), we can state:
\begin{theorem}
If M be a conformally flat perfect fluid K\"ahler space-time manifold satisfying the Einstein equation with cosmological constant then the manifold is infinitesimally spatially isotropic relative to the velocity vector field $\rho$.  
\end{theorem}

\section{Weakly symmetric perfect fluid K\"ahler space-time manifold}
If $M$ be a weakly symmetric K\"ahler manifold then we have
\begin{equation}
R(\overline Y,\overline Z,U,V)=R(Y,Z,U,V).
\end{equation}
Taking covariant derivative of (4.1), we get easily
\begin{equation}
(\nabla_XR)(\overline Y,\overline Z,U,V)=(\nabla_XR)(Y,Z,U,V).
\end{equation}
Using (1.3) in (4.2), we have
\begin{equation}
\begin{split}
\omega(Y)R(X,Z,U,V)+\omega(Z)R(Y,X,U,V)&=\omega(\overline Y)R(X,\overline Z,U,V)\\&+\omega(\overline Z)R(\overline Y,X,U,V).
\end{split}
\end{equation}
Putting $Z=U=e_i, 1\leq i\leq 4$ in (4.3) and taking summation over $i$, we obtained
\begin{equation}
\omega(Y)S(X,V)-R(Y,X,V,\rho)=\omega(\overline Y)S(X,\overline V)+R(\overline Y,X,V,\overline\rho).
\end{equation}
By using (2.8), equation (4.4) implies
\begin{equation}
\frac{r}{4}\omega(Y)g(X,V)-R(Y,X,V,\rho)=\frac{r}{4}\omega(\overline Y)g(X,\overline V)+R(\overline Y,X,V,\overline\rho).
\end{equation}
Putting $X=V=e_i, 1\leq i\leq 4$ in (4.5) and taking summation over $i$, we get
\begin{equation}
S(Y,\rho)=\frac{r}{2}\omega(Y).
\end{equation}
or
\begin{equation}
S(Y,\rho)=\frac{r}{2}g(Y,\rho).
\end{equation}
Replacing $\rho$ by $\overline\rho$ in (4.7), we can write
\begin{equation}
S(Y,\overline\rho)=\frac{r}{2}g(Y,\overline\rho).
\end{equation}
Hence from (4.7) and (4.8), we can state:
\begin{theorem}
If M be a weakly symmetric perfect fluid K\"ahler space-time manifold satisfying the Einstein equation with cosmological constant then $\rho$ and $\overline\rho$ are the eigen vector of the Ricci tensor $S$ with respect to eigen value $\frac{r}{2}$. 
\end{theorem}
Now putting $X=V=e_i, 1\leq i\leq 4$ in (1.3) and summing over $i$, we have
\begin{equation}
\begin{split}
(divR)(Y,Z)U=& R(Y,Z,U,\alpha)+\omega(Y)S(Z,U)\\&+\omega(Z)S(Y,U)+R(Y,Z,U,\rho).
\end{split}
\end{equation}
From (2.8), it can be easily obtained
\begin{equation}
(\nabla_YS)(Z,U)=\frac{r}{4}(\nabla_Yg)(Z,U)=0.
\end{equation}
 Using (4.10) in Bianchi second identity, we can write
\begin{equation}
(divR)(Y,Z)U=(\nabla_YS)(Z,U)-(\nabla_ZS)(Y,U)=0.
\end{equation}
By using (4.11), the equation (4.9) gives
\begin{equation}
R(Y,Z,U,\alpha)+\omega(Y)S(Z,U)-\omega(Z)S(Y,U)+R(Y,Z,U,\rho)=0.
\end{equation}
Putting $Z=U=e_i, 1\leq i\leq 4$ in (4.12) and taking summation over $i$, we have
\begin{equation}
S(Y,\alpha)+r\omega(Y)=0.
\end{equation}
Using (2.8) in (4.13), we can write
\begin{equation}
\frac{r}{4}g(Y,\alpha)+rg(Y,\rho)=0.
\end{equation}
Replacing Y by $\rho$ in (4.14), we get
\begin{equation}
r[g(\alpha,\rho)-4]=0,
\end{equation}
which implies either $r=0$ or $g(\alpha,\rho)=4$.\\
Hence, from above discussion we have
\begin{theorem}
If M be a weakly symmetric perfect fluid K\"ahler space-time manifold satisfying the Einstein equation with cosmological constant then either the manifold is of zero scalar curvature or the associated vector fields $\alpha$ and $\rho$ are related by $g(\alpha,\rho)=4$. 
\end{theorem}

\section{Weakly Ricci symmetric perfect fluid K\"ahler space-time manifold}

From equation (4.10) it is clear that if the manifold is perfect fluid K\"ahler space-time manifold then $(\nabla_XS)(Y,Z)=0$, therefore from (1.4), we have
\begin{equation}
A(X)S(Y,Z) +\omega(Y)S(Z,X)+\omega(Z)S(X,Y)=0,
\end{equation}
for weakly Ricci symmetric perfect fluid K\"ahler space-time manifold.
Using (2.8) in (5.1), we can write
\begin{equation}
\frac{r}{4}[A(X)g(Y,Z) +\omega(Y)g(Z,X)+\omega(Z)g(X,Y)]=0.
\end{equation}
equation (5.2) implies either scalar curvature $r=0$ or 
\begin{equation}
A(X)g(Y,Z) +\omega(Y)g(Z,X)+\omega(Z)g(X,Y)=0.
\end{equation}
Now, if (5.3) holds then by replacing Y and Z by $\overline Y$ and $\overline Z$ in (5.3) and using (1.6), we get
\begin{equation}
A(X)g(Y,Z) +\omega(\overline Y)g(\overline Z,X)+\omega(\overline Z)g(X,\overline Y)=0.
\end{equation}
Subtracting (5.3) from (5.4), we have
\begin{equation}
\omega(\overline Y)g(X,\overline Z)-\omega(Y)g(X,Z)+\omega(\overline Z)g(X,\overline Y)-\omega(Z)g(X,Y)=0.
\end{equation} 
Putting $X=Z=e_i, 1\leq i\leq 4$ in (5.5) and taking summation over $i$, we get
\begin{equation}
\omega(Y)=0,
\end{equation}
which is not possible because $g(\rho,\rho)=-1$.\\
Hence, we can state that:
\begin{theorem}
There does not exist a weakly Ricci symmetric perfect fluid K\"ahler space-time manifold satisfying the Einstein equation with cosmological constant having non-zero scalar curvature tensor.
\end{theorem}

{}

\end{document}